\documentclass[11pt]{article}
\usepackage{amssymb}
\usepackage{amsmath}

\newtheorem{theo}{Theorem}
\newtheorem{prop}{Proposition}
\newtheorem{lemm}{Lemma}
\newtheorem{coro}{Corollary}

\newcommand{\cqfd}
{%
\mbox{}%
\nolinebreak%
\hfill%
\rule{2mm}{2mm}%
\medbreak%
\par%
}
\newfont{\gothic}{eufb10}
\date{\empty}
\begin{document}
\title{On the homotopy types of K\"ahler manifolds and the birational Kodaira problem}
\author{Claire Voisin\\ Institut de math{\'e}matiques de Jussieu, CNRS,UMR
7586} \maketitle \setcounter{section}{-1}
\section{Introduction}
In small dimensions, it is known that K\"ahler compact manifolds are deformation equivalent
to smooth projective complex varieties. In dimension $2$, this follows
from the following theorem:
\begin{theo} (Kodaira \cite{kodairasurf}) Any compact K\"ahler surface admits
small deformations which are projective.
\end{theo}
The so-called Kodaira problem left open by this result asked whether more generally
any compact K\"ahker manifold can be deformed to a projective complex manifold.

Recently, we solved negatively this question by constructing,
in any dimension $n\geq4$,  examples of compact K\"ahler
manifolds, which  do not deform to projective complex manifolds,
as a consequence of the following stronger statement concerning
the topology of K\"ahler compact manifolds:

\begin{theo} (Voisin, \cite{voisin}) In any dimension $n\geq4$, there
are examples of compact K\"ahler
manifolds, which  do not  have the homotopy type of
projective complex manifolds.
\end{theo}

However, these examples were obtained starting either from certain complex tori or from
self-products of certain generalized Kummer varieties, and then blowing-up them along
adequate subsets.

Hence, all these examples are bimeromorphically equivalent to other
complex manifolds which satisfy the property of deforming
 to projective complex manifolds, namely complex tori, or self-products of
generalized Kummer
varieties.

The following question, which was asked to me by N. Buchdahl, F. Campana, S.-T. Yau,
and can be considered as a birational version of the Kodaira problem, is thus
quite natural:

\vspace{0,5cm}

{\bf Question.} {\it Let $X$ be a compact K\"ahler manifold. Does there exist a bimeromorphic model
$X'$ of $X$ which deforms to a projective complex manifold?}

\vspace{0,5cm}

In this paper, we show that the answer to this question is again no, which follows from
the following  stronger statement:
\begin{theo}\label{main} In any even dimension
$\geq 8$, there exist compact K\"ahler manifolds $X$, such that no compact
bimeromorphic model
$X'$ of $X$ has the homotopy type of a projective complex manifold.

\end{theo}
In this statement, we can in fact replace ``homotopy type'' with
``rational homotopy type'', that is ``rational cohomology ring''
(see Theorem \ref{maincoh}). Indeed, the whole discussion deals with
the Hodge structure on rational cohomology and the (non)-existence
of polarizations on them.

\vspace{0,5cm}

 {\bf Acknowledgements.} This work was started at Harvard
University; I thank S.-T. Yau for his invitation there and for his
interest in the work \cite{voisin}. I also thank N. Buchdahl, F.
Campana, S.-T. Yau for asking the question above.
\section{Construction of examples\label{sec1}}
We start as in \cite{voisin}, namely, we consider
$n$-dimensional complex tori $T$ admitting an endomorphism
$$\phi_T:T\rightarrow T$$
satisfying the following property (*). We can write $T$ as $\Gamma_{\mathbb{C}}/(\Gamma\oplus\Gamma')$,
where $\Gamma$ is a rank $2n$ lattice, $\Gamma_{\mathbb{C}}=\Gamma\otimes\mathbb{C}$ and
$\Gamma'$ is a complex subspace of $\Gamma_{\mathbb{C}}$ of rank $n$ such that
$$\Gamma'\oplus\overline{\Gamma}'=\Gamma_{\mathbb{C}}.$$
Let $\phi$ be the endomorphism $\phi_{T*}$ of $H_1(T,\mathbb{Z})=\Gamma$.
Clearly $\Gamma'$ has to be an eigenspace of $\phi_\mathbb{C}$, so no eigenvalue of $\phi$ can
be real.
The condition (*) is the following:

\vspace{0,5cm}

(*) {\it The characteristic polynomial of
$\phi$ (which has integer coefficients), has $2n$ distinct roots
(the eigenvalues of $\phi$) and its Galois group over $\mathbb{Q}$
 acts as the symmetric group of $2n$ letters on them.}

\vspace{0,5cm}

In the sequel, we will need to assume that the dimension $n$ of $T$ is at least $4$.
We make now the following construction. Let $\widehat{T}$ be the dual torus of $T$, namely
$$\widehat{T}=\Gamma_\mathbb{C}^*/(\Gamma^*\oplus\Gamma'^{\perp}).$$
Geometrically, $\widehat{T}$ is the torus
$$Pic^0(T)=H^1(T,\mathbb{C})/(H^{1,0}(T)\oplus H^1(T,\mathbb{Z}))$$
which is the group of topologically trivial holomorphic line bundles on $T$ up to isomorphism.

There exists on $T\times\widehat{T}$ the so-called Poincar\'e line bundle
$\mathcal{P}$ which is
uniquely characterized by the following properties :

- For any $t\in \widehat{T}$ parameterizing a line bundle $L_t$ on $T$, we have
$$L_t\cong\mathcal{P}_{\mid T\times t}.$$

- The restriction $\mathcal{P}_{\mid 0\times \widehat{T}}$ is trivial.

In fact $\mathcal{P}$ is constructed as follows : first of all, its first Chern class
$$c_1(\mathcal{P})\in NS(T\times\widehat{T}):=H^{1,1}(T\times \widehat{T})\cap
H^2(T\times \widehat{T},\mathbb{Z})$$ is the identity
$$id_{H^1(\widehat{T})}\in H^1(T,\mathbb{Z})\otimes
H^1(\widehat{T},\mathbb{Z})\subset H^2(T\times \widehat{T},\mathbb{Z}),$$
which is easily seen to be of Hodge type $(1,1)$.
Next the uniqueness of $\mathcal{P}$ is forced by the conditions
$$\mathcal{P}_{\mid 0\times \widehat{T}}\cong\mathcal{O}_{\widehat{T}},\,
\mathcal{P}_{\mid {T}\times0}\cong\mathcal{O}_{{T}}.$$ Next, because
$T$ admits the endomorphism $\phi_T$, we also have the line bundle
$$\mathcal{P}_{\phi}:=(\phi,Id)^*\mathcal{P}.$$

We now make the following construction:Over $T\times\widehat{T}$, consider the rank $2$
 vector bundles
$$E=\mathcal{P}\oplus\mathcal{P}^{-1},\,\,E_\phi=\mathcal{P}_\phi\oplus\mathcal{P}_\phi^{-1}$$
and the corresponding associated projective bundles $\mathbb{P}(E),\,\mathbb{P}(E_\phi)$.
The two commuting involutions
$(-Id,Id)$ and $(Id,-Id)$ of
$T\times\widehat{T}$ lift to commuting involutions $i,\,\hat{i}$, resp. $i_\phi,\,\hat{i}_\phi$
acting
on $E$ resp. $E_\phi$, since we have isomorphisms
$$(-Id,Id)^*\mathcal{P}\cong\mathcal{P}^{-1},\,(Id,-Id)^*\mathcal{P}\cong\mathcal{P}^{-1},$$
$$(-Id,Id)^*\mathcal{P}_\phi\cong\mathcal{P}_\phi^{-1},\,
(Id,-Id)^*\mathcal{P}_\phi\cong\mathcal{P}_\phi^{-1},$$
which can be made canonical by a choice of trivialization
$$\mathcal{P}_{\mid (0,0)}\cong \mathbb{C},$$
$(0,0)$ being a fixed point of both $(Id,-Id)$ and $(-Id,Id)$.

The compact K\"ahler manifold we shall consider is the following:We
start with the fibered product
$$\mathbb{P}(E)\times_{T\times\widehat{T}}\mathbb{P}(E_\phi).$$
It admits the commuting involutions
$$(i,i_\phi),\, (\hat{i},\hat{i}_\phi)$$ over $(-Id,Id)$, $(Id,-Id)$ respectively.
The quotient $Q$ of
$\mathbb{P}(E)\times_{T\times\widehat{T}}\mathbb{P}(E_\phi)$ by the
group $\mathbb{Z}/2\mathbb{Z}\times\mathbb{Z}/2\mathbb{Z}$ generated
by these involutions is singular along the non free locus of this
action, but the quotient admits a K\"ahler compact
desingularization. For example, one can start by desingularizing the
quotient
$\mathbb{P}(E)\times_{T\times\widehat{T}}\mathbb{P}(E_\phi)/(i,i_\phi)$
by blowing-up the fixed locus of $(i,i_\phi)$ and then taking the
quotient of the blown-up variety  by the natural involution which
lifts $(i,i_\phi)$. The result is smooth K\"ahler and by naturality
$(\hat{i},\hat{i}_\phi)$ acts on it as an involution. Then one can
desingularize in the same way the quotient of this new variety by
$(\hat{i},\hat{i}_\phi)$.

Our compact K\"ahler manifold  $X$ will be any K\"ahler
desingularization of this quotient.

Note that, if $K$ is the Kummer variety of $T$, namely the
desingularization of the quotient of $T$ by the $-Id$ involution,
obtained by blowing-up the images of the $2$-torsion points of $T$,
and similarly $\widehat{K}$ is the Kummer variety of $\widehat{T}$,
then over $K_0\times \widehat{K}_0$, $X$ is a
$\mathbb{P}^1\times\mathbb{P}^1$-bundle, where $K_0$ is the open set
$T_0/\pm Id$ of $K$, with
$$T_0:=T\smallsetminus 2-torsion\,\,points,$$
 and similarly for
$\widehat{K}_0$.

The next sections will be devoted to the proof of the following Theorem:
\begin{theo}\label{property} Let $X'$ be any compact complex manifold bimeromorphically
equivalent to
$X$. Then $X'$ does not have the homotopy type of a complex projective manifold.
\end{theo}

\section{Some results on the cohomology ring of $X'$}
We plan to show in fact the following slightly stronger
result:
\begin{theo}\label{maincoh} Let $X'$ be any compact complex manifold
bimeromorphically equivalent to $X$, and let $Y$ be a K\"ahler
compact manifold. Assume there is an isomorphism of graded algebras:
$$\gamma:H^*(Y,\mathbb{Q})\cong H^*(X',\mathbb{Q}).$$
Then $Y$ is not projective.
\end{theo}
In other words, Theorem \ref{property} is true for rational homotopy
type rather than homotopy type, since it is known that the rational
homotopy type of a compact K\"ahler manifold is determined by its
rational cohomology algebra (see \cite{DGM}).

This section will be devoted to the study of the cohomology ring
 of any compact complex manifold $X'$ given as in Theorem \ref{maincoh}.
  The proof of Theorem \ref{maincoh} will be given in the next section,
 following the same line as \cite{voisin}, section 3.

 Recall that $X$ admits a holomorphic map
 $$q:X\rightarrow (T/\pm Id)\times(\widehat{T}/\pm Id),$$
 obtained by composing the desingularization map
 $$ X\rightarrow \mathbb{P}(E)\times_{T\times\widehat{T}}\mathbb{P}(E_\phi)/
 <(i,i_\phi),(\hat{i},\hat{i}_\phi)>$$
 with the natural map
 $$\mathbb{P}(E)\times_{T\times\widehat{T}}\mathbb{P}(E_\phi)/<(i,i_\phi),(\hat{i},\hat{i}_\phi)>
 \rightarrow T\times\widehat{T}/<(-Id,Id),(Id,-Id)>.$$
For simplicity of notations,
we shall assume in the sequel that  our $X$ in section \ref{sec1}
has been chosen so that
$q$ extends to a holomorphic map
$$\overline{q}:X\rightarrow K\times\widehat{K},$$
which can always be achieved by a bimeromorphic transformation.
 \begin{lemm}\label{le1} Let $\psi:X'\dasharrow X$ be any
 bimeromorphic map. Then $q\circ \psi$ is holomorphic.
 \end{lemm}
 {\bf Proof.} The complex manifold $T\times \widehat{T}$ does not
 contain any closed complex curve. Indeed, it suffices to prove this for $T$ or
 $\widehat{T}$. Now, the cohomology class
 $[C]$ of such a curve would be a non zero Hodge class of degree $2n-2$
 on $T$, resp. $\widehat{T}$,
 or equivalently, a non-zero Hodge class in $H^2(\widehat{T},\mathbb{Q})$,
 resp. $H^2({T},\mathbb{Q})$.
 But in \cite{voisin}, Remark 3, we  proved that the existence of
 $\phi_T$, resp. ${\phi}_{\widehat{T}}$ prevents the existence of such a Hodge
 class.

 It follows that the quotient
 $(T/\pm Id)\times(\widehat{T}/\pm Id)$ does not contain any rational
 curve, and by desingularization of meromorphic maps with value in
 compact complex manifolds, this is enough
 to conclude that $q\circ\psi$ has to be holomorphic.
 \cqfd
It follows from Lemma \ref{le1} that $H^*(X',\mathbb{Q})$ contains a
subalgebra
$$A^*:=(q\circ\psi)^*H^*((T/\pm Id)\times(\widehat{T}/\pm
Id),\mathbb{Q})$$ which is isomorphic to $H^*((T/\pm Id)\times(\widehat{T}/\pm
Id),\mathbb{Q})$. Note that this last space is isomorphic to
$$H^*(T/\pm Id,\mathbb{Q})\otimes H^*(\widehat{T}/\pm
Id,\mathbb{Q}),$$ and that $H^*(T/\pm
Id,\mathbb{Q})=H^{even}(T,\mathbb{Q})$ (and similarly for
$\widehat{T}$).

We shall denote by
$A_1^*$, resp. $A_2^*$, the subalgebra
$(q\circ\psi)^* H^*(T/\pm Id\times 0,\mathbb{Q})$, resp.
$(q\circ\psi)^* H^*(0\times \widehat{T}/\pm Id ,\mathbb{Q})$.

Next, we note that the cohomology of $X'$ in degree $2$ is generated
over $\mathbb{Q}$ by $A^2$ and by degree $2$ Hodge classes. Indeed,
this is  true for $X$, because $X$ contains a Zariski open set which
is a $\mathbb{P}^1\times\mathbb{P}^1$-bundle over $K_0 \times
\widehat{K}_0$, and this implies easily that
$H^{2,0}(X)=H^0(X,\Omega_X^2)$ is equal to
$$\overline{q}^*H^{2,0}(K\times\widehat{K})=
{q}^*H^{2,0}((T/\pm Id)\times(\widehat{T}/\pm Id)).$$
 Next,
this property is invariant under meromorphic transformations,
hence if it is true for $X$, it is true
for $X'$.

Let now $D\subset H^2(X',\mathbb{Q})$ be the subspace generated by
degree $2$ Hodge classes. So we have
\begin{eqnarray}\label{decomp}H^2(X',\mathbb{Q})=D\oplus A^2, \end{eqnarray} because, by
\cite{voisin},  Remark 3,  we know that the presence of the endomorphism
$\phi_T$ of $T$ satisfying property (*) of section \ref{sec1}
implies that  $H^2(T,\mathbb{Q})$ has no non-zero Hodge class, and similarly for
$\widehat{T}$.
Furthermore, we have by definition
\begin{eqnarray}\label{kun}A^2=A_1^2\oplus A_2^2.
\end{eqnarray}
For $\alpha\in H^2(X',\mathbb{C})$, let
$$\alpha=\alpha_D+\alpha',\,\alpha'=\alpha_1+\alpha_2,$$ be
its decompositions given by (\ref{decomp}), (\ref{kun}).

A key role will be played by the following Proposition \ref{propcruciale}:

 We consider the algebraic subset
$Z\subset H^2(X',\mathbb{C})$
defined as
$$Z=\{\alpha\in H^2(X',\mathbb{C}),\,\,\alpha^2=0\,\,{\rm in}\,\, H^4(X',\mathbb{C})\}.$$
$Z$ contains the algebraic subsets $Z_1$, $Z_2$ defined as
$$ Z_1=\{\alpha\in H^2(X',\mathbb{C}),\,\alpha_2=0,\,\alpha_1^2=\alpha_1\alpha_D=\alpha_D^2=0
\,\,{\rm in}\,\, H^4(X',\mathbb{C})\},$$ resp.
$$ Z_2=\{\alpha\in H^2(X',\mathbb{C}),\,\alpha_1=0,\,\alpha_2^2=\alpha_2\alpha_D=0
=\alpha_D^2=0
\,\,{\rm in}\,\, H^4(X',\mathbb{C})\}.$$
\begin{prop}\label{propcruciale} Any irreducible component of
$Z_1$ (resp. $Z_2$) containing $$Z_{1,0}:=Z_1\cap
\{\alpha,\,\alpha_D=0\}$$ (resp. $Z_{2,0}:=Z_2\cap
\{\alpha,\,\alpha_D=0\}$) is an irreducible component of $Z$.

\end{prop}
{\bf Proof.}
The condition $\alpha^2=0$ writes as
\begin{eqnarray}\label{1eq}\alpha'^2+2\alpha_D\alpha'+\alpha_D^2=0.
\end{eqnarray}
Now we observe that $\alpha_D^2$ belongs to
$Hdg^4(X')\otimes\mathbb{C}$, where
$$Hdg^4(X'):=H^4(X',\mathbb{Q})\cap H^{2,2}(X').$$
Similarly, because the Hodge structure on $D$ is trivial, that is
purely of type $(1,1)$, $\alpha_D\alpha'$ belongs to
$N_2H^4(X')\otimes \mathbb{C}$, where $N_2H^4(X')$ is the maximal
rational sub-Hodge structure of $H^4(X',\mathbb{Q})$ which is of
Hodge level $2$.

(Here we recall that the level of a weight $k$ Hodge structure
$H,\,H_\mathbb{C}=\oplus_{p+q=k}H^{p,q}$, is the integer $$
Max\,\{p-q,\,H^{p,q}\not=0\}.$$ Thus a level $2$ sub-Hodge structure of a
weight $4$ Hodge structure, is a sub-Hodge structure
which has no $(4,0)$-term.)

 Equation (\ref{1eq}) thus implies
that $\alpha'^2$ belongs to $N_2A^4_\mathbb{Q}\otimes \mathbb{C}$,
where again $N_2$ means that we consider the maximal rational
sub-Hodge structure of level $2$.

Next we have the K\"unneth decomposition:
\begin{eqnarray}\label{decomp4}
A^4_\mathbb{Q}={A^4_1}_\mathbb{Q}\oplus {A_1^2}_\mathbb{Q}\otimes
{A_2^2}_\mathbb{Q}\oplus {A^4_2}_\mathbb{Q},
\end{eqnarray}
which is a decomposition into sub-Hodge structures of weight $4$.
We have the following:
\begin{lemm}\label{leirred4}${A^4_1}_\mathbb{Q}$ and $ {A^4_2}_\mathbb{Q}$ do not contain
non trivial sub-Hodge structure of Hodge level $2$.

\end{lemm}
{\bf Proof.} We use the fact that $n\geq 4$, so that
${A^4_1}_\mathbb{Q}$ and $ {A^4_2}_\mathbb{Q}$ are of Hodge level
$4$. Next we use the assumption (*) satisfied by $\phi$ to conclude
that $\phi_T^*$ acts in an irreducible way on
$\bigwedge^4H^1(T,\mathbb{Q})={A^4_1}_\mathbb{Q}$, and since the
action is via morphisms of Hodge structures, it must preserve
$N_2{A^4_1}_\mathbb{Q}$. Hence, because  $N_2{A^4_1}_\mathbb{Q}\not={A^4_1}_\mathbb{Q}$,
we conclude that
$N_2{A^4_1}_\mathbb{Q}=0$ and similarly $N_2{A^4_2}_\mathbb{Q}=0$.
\cqfd From the fact that
$\alpha'^2=\alpha_1^2+2\alpha_1\alpha_2+\alpha_2^2 \in
N_2A^4_\mathbb{Q}\otimes \mathbb{C}$, from the decomposition
(\ref{decomp4}) into sub-Hodge structures and from Lemma
\ref{leirred4}, we conclude that
$$\alpha_1^2=0,\,\alpha_2^2=0.$$
Thus our initial equation (\ref{1eq}) becomes
\begin{eqnarray}\label{2eq}2\alpha_1\alpha_2+2\alpha_D\alpha'+\alpha_D^2=0.
\end{eqnarray}
\cqfd
This equation implies as already noticed that
$\alpha_1\alpha_2$ belongs to the space
$$N_2({A_1^2}_\mathbb{Q}\otimes
{A_2^2}_\mathbb{Q})\otimes\mathbb{C}.$$ In fact one can say more:
indeed, note that the Hodge structure on $$D\cdot
A^2_\mathbb{Q}+D^2\subset H^4(X',\mathbb{Q})$$ is  the quotient of a
direct sum of Hodge structures of level $2$ isomorphic either to
$A^2_{1\mathbb{Q}}$ or to $A^2_{2\mathbb{Q}}$ or to  a trivial Hodge
structure.

Thus  condition (\ref{2eq}) implies that $\alpha_1\alpha_2$ has in
fact to belong to the space
$$N'_2({A_1^2}_\mathbb{Q}\otimes
{A_2^2}_\mathbb{Q})\otimes\mathbb{C},$$ where $N'_2$ means {\it the
maximal sub-Hodge structure of level $2$, which is a subquotient of
a sum of copies of
 $A^2_{1\mathbb{Q}}$ or  $A^2_{2\mathbb{Q}}$ or a trivial Hodge structure}.

On the other hand, the Hodge structures on  $A^2_{1\mathbb{Q}}$ or
$A^2_{2\mathbb{Q}}$  are simple, that is do not contain any
non-trivial sub-Hodge structure. To see this last point, assume that
there is a proper non-zero simple sub-Hodge structure
$$H\subset H^2(T,\mathbb{Q}).$$
As the endomorphism of Hodge structure $\phi_T^*$ acts transitively
on $H^2(T,\mathbb{Q})$, it follows that $H^2(T,\mathbb{Q})$ must
then be isomorphic to a sum of copies of $H$,
$$\exists k>1,\,\,H^2(T,\mathbb{Q})\cong H^k.$$
But then $H^2(T,\mathbb{Q})$ admits a projector which is an
endomorphism of Hodge structure. This contradicts the fact, noted at
the end of the proof of Lemma \ref{plustard2}, that the algebra
$End\,H^2(T,\mathbb{Q})$ is generated by $\phi_T^*$, and thus does
not contain projectors by condition (*).

 Note also that  the Hodge structures on  $A^2_{1\mathbb{Q}}$ and
$A^2_{2\mathbb{Q}}$ are not isomorphic, as shown by Lemma
\ref{plustard1} below.

Thus it follows that $N'_2({A_1^2}_\mathbb{Q}\otimes
{A_2^2}_\mathbb{Q})$  is in fact equal to   the maximal sub-Hodge
structure of level $2$ of ${A_1^2}_\mathbb{Q}\otimes
{A_2^2}_\mathbb{Q}$, which is  a sum of copies of
 $A^2_{1\mathbb{Q}}$ or  $A^2_{2\mathbb{Q}}$ or a trivial Hodge structure.

 We have  the following Lemma:
\begin{lemm} \label{plustard1} There are no non zero morphism of
Hodge structures (of bidegree $(1,1)$) from $A^2_{1\mathbb{Q}}$ or
$A^2_{2\mathbb{Q}}$ to ${A_1^2}_\mathbb{Q}\otimes
{A_2^2}_\mathbb{Q}$.
\end{lemm}
Admitting this Lemma, we conclude that in fact $\alpha_1\alpha_2$
has to belong to $$Hdg({A_1^2}_\mathbb{Q}\otimes
{A_2^2}_\mathbb{Q})\otimes \mathbb{C},$$ where $Hdg$ means {\it the
subspace of rational Hodge classes}. We have next the following
Lemma, the proof of which we shall also postpone:
\begin{lemm} \label{plustard2}There are (up to a coefficient) finitely many elements
$$\beta\in Hdg({A_1^2}_\mathbb{Q}\otimes
{A_2^2}_\mathbb{Q})\otimes \mathbb{C}$$ which are of rank $1$, that
is of the form $\alpha_1\alpha_2$ as above.
\end{lemm}
We then conclude as follows: from the above analysis, we conclude that
for $\alpha\in Z$, we have either $\alpha_1\not=0,\,\alpha_2\not=0$ and then
$\alpha_1\alpha_2$ has to be proportional to one of the finitely many
$\beta$ of Lemma \ref{plustard2}, or
one of $\alpha_1$, or $\alpha_2$ has to be $0$.

We claim that in this last case, $\alpha$ belongs to $Z_2$ or $Z_1$
respectively. Indeed, we know that in any case
$$\alpha_1^2=0,\,\alpha_2^2=0.$$
Assume $\alpha_2=0$. Equation (\ref{2eq}) thus becomes:
$$2\alpha_D\alpha_1+\alpha_D^2=0.
$$
But this implies that
$$\alpha_1\alpha_D=0,\,\alpha_D^2=0.$$
Indeed, $\alpha_D^2$ belongs to
$Hdg^4(X')\otimes\mathbb{C}$ while
$\alpha_1\alpha_D$ belongs to the space
$$N''_2H^4(X',\mathbb{Q})\otimes\mathbb{C},$$
defined as the maximal sub-Hodge structure of $H^4(X',\mathbb{Q})$
isomorphic to a subquotient of some power of $A^2_{1\mathbb{Q}}$. By
the same simplicity argument as before,
$N''_2H^4(X',\mathbb{Q})\otimes\mathbb{C}$ is also the maximal
sub-Hodge structure of $H^4(X',\mathbb{Q})$ isomorphic to  some
power of $A^2_{1\mathbb{Q}}$.But the intersection
$$Hdg^4(X')\cap N''_2H^4(X',\mathbb{Q})$$
has to be zero, since there is no non zero Hodge class in
$A^2_{1\mathbb{Q}}$. Thus also
$$Hdg^4(X')\otimes\mathbb{C}\cap N''_2H^4(X',\mathbb{Q})\otimes\mathbb{C}$$
has also to be $0$. Hence we proved that
$2\alpha_D\alpha_1+\alpha_D^2=0$ implies that
$\alpha_1\alpha_D=0,\,\alpha_D^2=0$.

In conclusion, we proved that
$Z$ is the set-theoretic  union of $Z_1,\,Z_2$, and of a set which projects to a finite set
of lines in $A^2_{1\mathbb{C}}$ and $A^2_{2\mathbb{C}}$.

 Let now $Z'_1$  be an irreducible component of $Z_1$ which contains
$Z_{1,0}$. Suppose it is not an irreducible component of $Z$. This
means that there exists an irreducible component $Z'$ of $Z$
containing $Z_{1}'$, not contained in $Z_1$, such that
$Z'\smallsetminus Z'_{1}$ is dense in $Z'$. So $Z'\smallsetminus
Z'_{1}$ has to project in a dominant way onto $A^2_{1\mathbb{C}}$,
which contradicts the fact that $Z'\smallsetminus Z'_{1}$ has to be
contained in the union of $Z_2$, which projects to $0$ in
$A^2_{1\mathbb{C}}$, and of a set which projects to a finite union
of lines  in $A^2_{1\mathbb{C}}$. Thus Proposition
\ref{propcruciale} is proved, assuming Lemmas \ref{plustard1} and
\ref{plustard2}. \cqfd {\bf Proof of Lemma \ref{plustard1}.} Recall
that
$$A_{1\mathbb{Q}}^2=\bigwedge^2H^1(T,\mathbb{Q})\cong\bigwedge^2\Gamma^*_\mathbb{Q},$$
$$A_{2\mathbb{Q}}^2=\bigwedge^2H^1(\widehat{T},\mathbb{Q})\cong\bigwedge^2\Gamma_\mathbb{Q}.$$
We have the endomorphisms $\phi_T,\,\phi_{\widehat{T}}$ acting
respectively on the complex tori $T$ and $\widehat{T}$, and the
induced action $\phi_T^*,\,\phi_{\widehat{T}}^*$ on
$H^2(T,\mathbb{Q})$, resp. $H^2(\widehat{T},\mathbb{Q})$, identify
to $\wedge^2 {^t\phi} ,\,\wedge^2 \phi$ respectively.

Let $\lambda_1,\ldots,\lambda_{2n}$ be the $2n$-eigenvalues of
$\phi$ on $\Gamma_\mathbb{C}$. Let $e_1,\ldots,e_{2n}$ be a
corresponding basis of eigenvectors of $\Gamma_\mathbb{C}$, and let
$e_i^*$ be the dual basis of $\Gamma_\mathbb{C}^*$. We choose the
ordering in such a way that $\Gamma'$ (see section \ref{sec1}) is
generated by $e_1,\ldots,e_n$. In other words, $e_i\in
H^1(\widehat{T},\mathbb{C})$ have Hodge type $(1,0)$ for $i\leq n$
and $e_i^*\in H^1(T,\mathbb{C})$ have Hodge type $(1,0)$ for $i>n$.

We want to study the Hodge classes in
$$A_{1,\mathbb{Q}}^{2*}\otimes A_{1,\mathbb{Q}}^{2}\otimes A_{2,\mathbb{Q}}^{2}$$
$$= \bigwedge^2\Gamma_\mathbb{Q}\otimes \bigwedge^2\Gamma_\mathbb{Q}^*
\otimes\bigwedge^2\Gamma_\mathbb{Q},$$ which we consider as a weight
$6$ Hodge structure, so the classes we search are the rational
classes of Hodge type $(3,3)$.

This space
$$S:=Hdg(\bigwedge^2\Gamma_\mathbb{Q}\otimes \bigwedge^2\Gamma_\mathbb{Q}^*
\otimes\bigwedge^2\Gamma_\mathbb{Q})$$is stable under the action of
the three commuting morphisms of Hodge structures
$$\wedge^2\phi\otimes Id\otimes Id,\, Id\otimes\wedge^2 {^t\phi}\otimes Id,\,
Id\otimes Id\otimes\wedge^2\phi.$$ It follows that the complexified
space $S_\mathbb{C}$ is generated by eigenvectors for these actions,
namely elements of the form
\begin{eqnarray} e_i\wedge e_j\otimes e_k^*\wedge e_l^*\otimes e_r\wedge e_s
\label{element1}.
\end{eqnarray}
For $a,\,b\,c\in \mathbb{Z}$, consider the endomorphism
$$\Phi_{abc}:=\wedge^2\phi\otimes\wedge^2
{^t\phi}\otimes\wedge^2\phi$$ of
$\bigwedge^2\Gamma_\mathbb{Q}\otimes \bigwedge^2\Gamma_\mathbb{Q}^*
\otimes\bigwedge^2\Gamma_\mathbb{Q} $. $\Phi_{abc}$ is diagonal in
the basis given by the elements (\ref{element1}), with corresponding
eigenvalues
$$(\lambda_i\lambda_j)^a(\lambda_k\lambda_l)^b(\lambda_r\lambda_s)^c.$$
 The Galois group of the field
$K=\mathbb{Q}[\lambda_1,\ldots,\lambda_{2n}]$ over $\mathbb{Q}$ acts
on the $\lambda_i$ and has to leave stable the set $E_{abc}$ of
eigenvalues of $\Phi_{abc}$ on $S$, since $S$ is defined over
$\mathbb{Q}$. On the other hand,  we know that this Galois group is
 the symmetric group $\mathfrak{S}_{2n}$ on $2n$ letters acting on the
 $\lambda_i$'s. Thus we conclude that if
 $$(\lambda_i\lambda_j)^a(\lambda_k\lambda_l)^b(\lambda_r\lambda_s)^c\in
 E_{abc},$$
 then also
 $$(\lambda_{\sigma(i)}\lambda_{\sigma(j)})^a(\lambda_{\sigma(k)}
 \lambda_{\sigma(l)})^b(\lambda_{\sigma(r)}\lambda_{\sigma(s)})^c
 \in
 E_{abc}.$$
 But for an adequate choice of $a,\,b,\,c$ the map
 $$(\{i,j\},\{k,l\},\{r,s\})\mapsto
 (\lambda_i\lambda_j)^a(\lambda_k\lambda_l)^b(\lambda_r\lambda_s)^c$$
 is injective.
 Thus we conclude that if (\ref{element1}) belongs
 to $S_\mathbb{C}$, so does
\begin{eqnarray} e_{\sigma(i)}\wedge e_{\sigma(j)}\otimes e_{\sigma(k)}^*\wedge
 e_{\sigma(l)}^*\otimes
e_{\sigma(r)}\wedge e_{\sigma(s)} \label{element2}.
\end{eqnarray}
As $S_\mathbb{C}$ is contained in the $(3,3)$-part of
$$A_{1,\mathbb{C}}^{2*}\otimes A_{1,\mathbb{C}}^{2}\otimes A_{2,\mathbb{C}}^{2}
$$
we see that (\ref{element2}) has to be of Hodge type $(3,3)$ for any
permutation $\sigma\in\mathfrak{S}_{2n}$.

 But as $n\geq4$, it is immediate that we can
always find $\sigma$ in such a way that (\ref{element2}) has Hodge
type $(4,2)$ (eg choose $i,\,j,\,r,\,s$ in $\{1,\ldots ,n\}$).

Thus an element (\ref{element1}) in
$$Hdg(\bigwedge^2\Gamma_\mathbb{Q}\otimes \bigwedge^2\Gamma_\mathbb{Q}^*
\otimes\bigwedge^2\Gamma_\mathbb{Q}) \otimes\mathbb{C}$$ does not
exist, which proves the lemma. \cqfd {\bf Proof of Lemma
\ref{plustard2}.} We study the space
$$Hdg({A_1^2}_\mathbb{Q}\otimes
{A_2^2}_\mathbb{Q})\otimes \mathbb{C}$$ exactly as in the previous proof.
The space ${A_1^2}_\mathbb{Q}\otimes
{A_2^2}_\mathbb{Q}$ identifies to
$$\bigwedge^2\Gamma^*_\mathbb{Q}\otimes\bigwedge^2\Gamma_\mathbb{Q},$$
as Hodge structures, where the $e_i^*\in \Gamma_\mathbb{C}^*$ have Hodge type
$(1,0)$ for $i> n$, while the $e_i\in \Gamma_\mathbb{C}$ have Hodge type
$(1,0)$ for $i\leq n$.

Again, the space $S:=Hdg({A_1^2}_\mathbb{Q}\otimes
{A_2^2}_\mathbb{Q})\otimes \mathbb{C}$, being stable under
$\wedge^2{^t\phi}\otimes Id$ and $Id\otimes \wedge^2\phi$ has to be
generated by eigenvectors for both of these commuting endomorphisms,
that is elements of the form :
$$ e_i^*\wedge e_j^*\otimes e_k\wedge e_l.$$
Because this space is defined over $\mathbb{Q}$, we conclude as in
the previous proof that it has to be stable under the action of
$\mathfrak{S}_{2n}$, which means that for any permutation $\sigma$
of $1,\ldots,2n$,
$$ e_{\sigma(i)}^*\wedge e_{\sigma(j)}^*\otimes e_{\sigma(k)}\wedge e_{\sigma(l)}$$
has to be of type $(2,2)$. Now this implies that up to permuting $i$
and $j$, one must have $i=k,\,j=l$. Indeed, if the four indices are
distinct, by changing  them by some $\sigma\in \mathfrak{S}_{2n}$,
we may arrange that $ e_{\sigma(i)}^*\wedge e_{\sigma(j)}^*\otimes
e_{\sigma(k)}\wedge e_{\sigma(l)}$ has Hodge type $(4,0)$, and if eg
$i=k$ but $j\not=l$, by changing  them by some $\sigma\in
\mathfrak{S}_{2n}$, we may arrange that $ e_{\sigma(i)}^*\wedge
e_{\sigma(j)}^*\otimes e_{\sigma(i)}\wedge e_{\sigma(l)}$ has Hodge
type $(3,1)$. Hence we have proved that
$Hdg({A_1^2}_\mathbb{Q}\otimes {A_2^2}_\mathbb{Q})\otimes
\mathbb{C}$ is generated by the elements $ e_i^*\wedge e_j^*\otimes
e_i\wedge e_j$ (in fact it has to be equal to the space generated by
these elements, which is nothing but the algebra generated over
$\mathbb{C}$ by $\phi_{\widehat{T}}^*$).

It is then clear that it contains (up to a scalar)
only finitely elements of rank $1$, namely the elements above.
\cqfd
Proposition \ref{propcruciale} is now fully proved.
Our next technical Lemma will be the following:
\begin{lemm}\label{produitdiv} Let $D_1\subset D$ be defined as
$$D_1=\{\alpha_D\in D,\,\alpha_D\alpha_1=0\,\,{\rm in}\,\,H^4(X',\mathbb{Q}),\,
\forall \alpha_1\in A^2_{1\mathbb{Q}}\}.$$
Then, if $\alpha_D\in D\otimes\mathbb{C}$ satisfies
$\alpha_D\alpha_1=0$ for one non zero $\alpha_1\in A^2_{1\mathbb{C}}$,
one has $\alpha_D\in D_1\otimes\mathbb{C}$.
\end{lemm}
{\bf Proof.} First of all, note that if
$$\psi':X''\rightarrow X'$$
is a proper surjective holomorphic map of degree $1$, with $X''$
smooth, and the result is true for $X''$, with $D$ replaced by the
space $Hdg^2(X'')$ and $A^2_{1\mathbb{Q}}$ by
$\psi'^*A^2_{1\mathbb{Q}}$, then it is also true for $X'$.

Indeed, such a map $\psi'$ induces an injective map $\psi'^*$ of
cohomology algebras, which sends $D$ in the space of Hodge classes
of degree $2$ on $X''$.

Recall now that $X'$ is bimeromorphic to a quotient of the
$\mathbb{P}^1\times\mathbb{P}^1$-bundle
$\mathbb{P}(E)\times_{T\times\widehat{T}}\mathbb{P}(E_\phi)$ over
$T\times\widehat{T}$. Hence there is a dominant meromorphic map
from $\mathbb{P}(E)\times_{T\times\widehat{T}}\mathbb{P}(E_\phi)$
to $X'$.

Using Hironaka's desingularization theorem and the previous
observation, we can thus reduce to the case where $X'$ is deduced from
 $W:=\mathbb{P}(E)\times_{T\times\widehat{T}}\mathbb{P}(E_\phi)$  by a
sequence of blow-ups.

We first prove that the result is true for $W$. The cohomology of
degree $2$ of $$W\stackrel{\tilde{q}}{\rightarrow} T\times
\widehat{T}$$ is a free module over the cohomology of $T\times
\widehat{T}$ generated by
$H^*(\mathbb{P}^1\times\mathbb{P}^1,\mathbb{Q})$. The space of
degree $2$ Hodge  classes $D$ on $W$ is the sum of two spaces, namely $D_0$
which has rank $2$ and is isomorphic by restriction to
$H^2(\mathbb{P}^1\times\mathbb{P}^1,\mathbb{Q})$ and $D_1$ which is
isomorphic via $\tilde{q}^*$ to the set of degree $2$ Hodge classes in $H^2(T\times
\widehat{T},\mathbb{Q})$. But we have by K\"unneth decomposition:
$$
H^2(T\times
\widehat{T},\mathbb{Q})=H^2(K\times\widehat{K},\mathbb{Q})\oplus
H^1(T,\mathbb{Q})\otimes H^1(\widehat{T},\mathbb{Q}),$$ and $D_1$ is
contained in the last factor.  One checks that
$D_1$ is generated over $\mathbb{Q}$
by $p:=c_1(\mathcal{P})$ and its pull-backs under
$(\phi_T^l)^*\otimes Id$.
We conclude from this that the product
map

$$D\otimes \tilde{q}^*(H^2(T,\mathbb{Q}),\mathbb{Q})\rightarrow
H^4(W,\mathbb{Q})$$ is injective, so that there is in fact nothing
to prove for $W$.

It remains now only to prove that if the statement  is true for $W$,
it is true for any complex manifold obtained by successive blow-ups
of $W$ along smooth centers. This is proved by induction on the
number of blow-ups. Assume it is true for $W_i$ and let
$\tau:W_{i+1}\rightarrow W_i$ be the blow-up of a  smooth
irreducible center $Z\subset W_i$. Then the set of degree $2$ Hodge classes
$D_{i+1}$ on $W_{i+1}$ is generated by $\tau^*D_i$ and the class
$e_Z$ of the exceptional divisor $E_Z$. Now, the study of the
cohomology ring of $W_{i+1}$ (see \cite{voisinbook} I, 7.3.3) shows that
if there is an equality
$$ e_Z\tau^*\alpha=0,\,modulo\,\,\tau^* H^*(W_i,\mathbb{C})$$
then in fact $ e_Z\tau^*\alpha=0$ in $H^*(W_{i+1},\mathbb{C})$.

Now suppose there is a relation $\alpha_D\alpha=0 $ in
$H^*(W_{i+1},\mathbb{C})$, where $\alpha_D\in
D_{i+1}\otimes\mathbb{C}$ and $\alpha\in
q^*H^2(T\times0,\mathbb{C})$. Writing
$$\alpha_D=\mu e_Z+\alpha'_D,$$
where $\mu\in \mathbb{C}$ and $\alpha'_D\in
\tau^*D_i\otimes\mathbb{C}$, we conclude using the previous remark
that
$$\mu e_Z\alpha=0\,\,{\rm in}\,\, H^*(W_{i+1},\mathbb{C}),$$
that is, either $\mu=0$, in which case we can apply the result for
$W_i$, or
$$e_Z\alpha=0\,\,{\rm in}\,\, H^*(W_{i+1},\mathbb{C}).$$
Since multiplication by the Hodge class $e_Z$ is a morphism of Hodge
structures from $H^2(T,\mathbb{Q})$ to $H^4(W_{i+1},\mathbb{Q})$,
its kernel is a sub-Hodge structure of $H^2(T,\mathbb{Q})$. So this
map is either injective or $0$, since the Hodge structure on
$H^2(T,\mathbb{Q})$ is simple, as already noticed before.

The conclusion is that, if there is one non-zero $\alpha$ satisfying $\alpha_D\alpha=0
$ in $H^*(W_{i+1},\mathbb{C})$ with  a coefficient $\mu\not=0$, we
find that $e_Z\alpha'=0$ in $H^*(W_{i+1},\mathbb{C})$, for any
$\alpha'\in q^* H^2(T\times0,\mathbb{Q})$, and that furthermore the
equality $\alpha_D\alpha=0 $ reduces to the equality $
\alpha'_D\alpha=0$, which holds already in $H^*(W_i,\mathbb{C})$.
Hence the result is proved by induction.

\cqfd
We will need also the following result.
\begin{lemm}\label{dercoh}
a) For any $d\in D\otimes \mathbb{C},\,\beta\in A^{4n-2}_\mathbb{C}\subset H^{4n-2}(X',\mathbb{C})$,
one has
$$d^3\beta=0\,\,{\rm in}\,\,H^{4n+4}(X',\mathbb{C})=\mathbb{C}.$$
b) The complex subspace $D\otimes \mathbb{C}\subset H^{2}(X',\mathbb{C})$ is an irreducible
component of the algebraic set
\begin{eqnarray}
\label{Z'}Z'=\{d\in H^2(X',\mathbb{C}),\,d^3\beta=0,\,\forall \beta\in A^{4n-2}_\mathbb{C}\}.
\end{eqnarray}
\end{lemm}
{\bf Proof.} $D$ is made of Hodge classes. So for any $d\in D$, the map
$$\alpha\mapsto d^3\alpha\in H^{4n+4}(X',\mathbb{Q})=\mathbb{Q}$$
is a Hodge class in $(A^{4n-2}_\mathbb{Q})^*=A^2_\mathbb{Q}$. But
 we already know that $A^2_\mathbb{Q}$ has no non zero Hodge classes.
This proves a).

Let $Z'_1\subset H^2(X',\mathbb{C})$ be an irreducible component of
the algebraic subset $Z'$ of (\ref{Z'}) containing strictly
$D\otimes \mathbb{C}$. Choose any point $d\in D\otimes \mathbb{C}$
and let
$$D'_\mathbb{C}:=T_{Z',d}\subset H^2(X',\mathbb{C}).$$
Since
$$D\oplus A^2_\mathbb{Q}=H^2(X',\mathbb{Q}),$$
and $D\otimes\mathbb{C}\subset D'_\mathbb{C}$, where the inclusion
is strict, there must be a non-zero element depending on $d$
\begin{eqnarray}\label{30sepalpha}
\alpha_d\in D'_\mathbb{C}\cap A^2_\mathbb{C}.
\end{eqnarray}
 This $\alpha_d$
satisfies the property that for any $\beta\in A^{4n-2}_\mathbb{C}$,
one has
\begin{eqnarray}\label{30sep}
d^2\alpha_d\beta=0\,\,{\rm in}\,\, H^{4n+4}(X',\mathbb{C}).
\end{eqnarray} We get a contradiction as follows: since $X'$ is in the
class $\mathcal{C}$, that is bimeromorphically equivalent to a
K\"ahler compact manifold, and the map
$$q\circ\psi:X'\rightarrow (T/\pm Id)\times(\widehat{T}/\pm Id)$$
is dominating with $4$-dimensional fiber, there is a $\mu\in H^2(X',\mathbb{C})$
such that
$$(q\circ\psi)_*\mu^2\not=0\,\,{\rm in}\,\,
H^{0}((T/\pm Id)\times(\widehat{T}/\pm Id),\mathbb{C})
\cong \mathbb{C}.$$ (Here we should work with $K\times \widehat{K}$
and desingularize  the map
$$q\circ\psi:X'\dasharrow K\times \widehat{K}$$
to be more rigorous on the definition of $(q\circ\psi)_*$.)

Now, write $\mu=d_1+\mu'$, with $d_1\in D\otimes\mathbb{C}$ and
$\mu'\in A^2_\mathbb{C}$. Then $$\mu^2=\mu'^2+2d_1\mu'+d_1^2,$$ so
that for any $\alpha\in A^2,\,\beta\in A^{4n-2}$,
$$\mu^2\alpha\beta=(\mu'^2+2d_1\mu'+d_1^2)\alpha\beta=d_1^2\alpha\beta.$$

Choose for $d$ the element $d_1$ above, and introduce $\alpha_{d_1}$
as in (\ref{30sepalpha}). Now, because $H^2((T/\pm
Id)\times(\widehat{T}/\pm Id),\mathbb{C})$ and $H^{4n-2}((T/\pm
Id)\times(\widehat{T}/\pm Id),\mathbb{C})$ are dual via the
cup-product and the isomorphism
$$H^{4n}((T/\pm Id)\times(\widehat{T}/\pm Id),\mathbb{C})=\mathbb{C},$$
there exists a $\beta\in H^{4n-2}((T/\pm Id)\times(\widehat{T}/\pm Id),\mathbb{C})$
such that
$$\alpha_{d_1}\beta\not=0 \,\,
{\rm in}\,\, H^{4n}((T/\pm Id)\times(\widehat{T}/\pm Id),\mathbb{C}).$$
Thus
$$\mu^2\alpha_{d_1}\beta\not=0\,\,{\rm in}\,\, H^{4n+4}(X',\mathbb{C}).$$
But we have just seen that
$$\mu^2\alpha_{d_1}\beta=d_1^2\alpha_{d_1}\beta.$$
The left hand side is non zero, while  the right hand side vanishes
by (\ref{30sep}), which proves b) by contradiction. \cqfd

We conclude this section with the proof of a Proposition concerning
the geometry of the bimeromorphic map $\psi:X'\dasharrow X$ which
will be essential in the sequel. Recall that we proved that the
meromorphic map
$$q\circ\psi:X'\dasharrow (T/\pm Id)\times(\widehat{T}/\pm Id)$$
is in fact holomorphic.
Let $X'_0:=(q\circ\psi)^{-1}(K_0\times\widehat{K}_0)$.
\begin{prop}\label{propess} There exists a dense Zariski open set
$U\subset K_0\times\widehat{K}_0$ such that denoting
$$X'_U:=(q\circ\psi)^{-1}(U),\,\,
X_U:=q^{-1}(U) , $$ the induced  meromorphic map
$$\psi:X'_U\dasharrow X_U$$
is holomorphic.
\end{prop}

In order to prove this proposition, we need to establish a few Lemmas saying
that $T\times \widehat{T}$ and $\mathbb{P}(E)\times_{T\times\widehat{T}}\mathbb{P}(E_\phi)$
contain very few closed analytic subsets.
They will be needed also later on in section \ref{sec3}.

\begin{lemm} \label{sub1}The only closed irreducible positive
dimensional proper analytic subsets of
$T\times\widehat{T}$ are of the
form $x\times\widehat{T}$, $x\in T$, or $T\times y$, $y\in \widehat{T}$.
\end{lemm}
{\bf Proof.} Indeed, note first that $T$ and $\widehat{T}$ do not
contain positive dimensional proper analytic subsets. This is
because they   both are simple tori which are not projective (see
\cite{voisin}), as guaranteed by the existence of $\phi_T$ and
$\phi_{\widehat{T}}$.

It follows that if $Z\subset T\times\widehat{T}$ is positive
dimensional proper irreducible and not of the above form, then it
must be \'etale over both $T$ and $\widehat{T}$ which implies that
the rational Hodge structures on $H^1(T,\mathbb{Q})$ and
$H^1(\widehat{T},\mathbb{Q})$ are isomorphic. But this is not the
case, as a consequence of Lemma \ref{plustard1}. \cqfd
\begin{lemm}\label{sub2} The only irreducible proper closed analytic subsets
of $\mathbb{P}(E)$ which dominate
$T\times\widehat{T}$ are the images $\Sigma_1,\,\Sigma_2$ of the two natural sections
$\sigma_1,\,
\sigma_2$ of $\mathbb{P}(E)$ corresponding to the splitting
$$E=\mathcal{P}\oplus \mathcal{P}^{-1},$$
and similarly for $\mathbb{P}(E_\phi)$.

\end{lemm}

{\bf Proof.} Indeed, let $Z\subset\mathbb{P}(E)$ be an hypersurface dominating
$T\times \widehat{T}$.
 Let us denote by $e:Z\rightarrow T\times \widehat{T}$
the generically finite map. Note that because of the description
above of the proper analytic subsets of $T\times \widehat{T}$, $Z$
has to contain a dense Zariski open set $Z_0$ which is an \'etale
cover of a Zariski open set $U\subset T\times \widehat{T}$, where
the complementary set of $U$ is an union of analytic subsets of the
form $x\times\widehat{T}$ or $T\times y$.

Next $Z$ induces a section of the induced $\mathbb{P}^1$-bundle
$\mathbb{P}(E)_Z:=e^*\mathbb{P}(E)$. Such a section is given by a
line bundle $\mathcal{L}$ over $Z$ and a surjective map
$$E^*=e^*\mathcal{P}\oplus e^*\mathcal{P}^{-1}\rightarrow \mathcal{L}.$$
If one of the two induced maps
$$e^*\mathcal{P}\rightarrow \mathcal{L}$$
or
$$e^*\mathcal{P}^{-1}\rightarrow \mathcal{L}$$
is zero, then $Z$ has to be contained in $\Sigma_1$ or $\Sigma_2$.
Otherwise, we find that both $e^*\mathcal{P}^{-1}\otimes
\mathcal{L}$ and $e^*\mathcal{P}\otimes \mathcal{L}$ have non-zero
sections. Note that, because $Z_0$ is an \'etale cover of an open
set of $T\times\widehat{T}$ whose complementary set has codimension
$\geq2$,  some power $\mathcal{L}^{\otimes k},\,k>0$ is equal to
$e^*(\mathcal{K})$ on $Z_0$, for some line bundle $\mathcal{K}$ on
$T\times\widehat{T} $. Furthermore,
$$e^*\mathcal{P}^{-k}\otimes \mathcal{L}^{\otimes k}=
e^*(\mathcal{P}^{-k}\otimes \mathcal{K})$$
and
$$e^*\mathcal{P}^{\otimes k}\otimes \mathcal{L}^{\otimes k}=e^*(
\mathcal{P}^{\otimes k}\otimes\mathcal{K})$$
 have non-zero sections
on $Z_0$.
It then follows that for some $m>0$, there are non-zero sections
of
$$\mathcal{P}^{-km}\otimes \mathcal{L}^{\otimes km}=\mathcal{P}^{-km}\otimes \mathcal{K}
^{\otimes m},$$
$$\mathcal{P}^{\otimes km}\otimes \mathcal{L}^{\otimes km}=\mathcal{P}^{\otimes km}
\otimes\mathcal{K}^{\otimes m},$$ on the open set $U$, hence on
$T\times\widehat{T}$ itself. But since $T\times\widehat{T}$  does
not contain hypersurfaces, these sections do not vanish anywhere,
from which one concludes that $\mathcal{P}^{-km}$ is isomorphic to
$\mathcal{P}^{km}$, which is not true since there cohomology classes
are different. This proves the Lemma for $\mathbb{P}(E)$ and the
result for $\mathbb{P}(E_\phi)$ follows. \cqfd
\begin{coro} \label{sub3}a) The only irreducible codimension $1$ analytic subsets of
$$\mathbb{P}(E)\times_{T\times\widehat{T}}\mathbb{P}(E_\phi)$$
which dominate $T\times\widehat{T}$ are of the form
$pr_1^{-1}\Sigma_i,\,i=1,\,2$ or $pr_2^{-1}\Sigma_i^\phi,\,i=1,\,2$.

b) The only irreducible codimension $2$ analytic subsets of
$$\mathbb{P}(E)\times_{T\times\widehat{T}}\mathbb{P}(E_\phi)$$
which dominate $T\times\widehat{T}$ are complete intersections
$$pr_1^{-1}\Sigma_i\cap pr_2^{-1}\Sigma_j^\phi,\,i=1,\,2,\,j=1,\,2.$$

\end{coro}
{\bf Proof.} Let $\mathcal{L}:=\mathcal{O}(Z)$, and let $
H=pr_1^*(\mathcal{O}_{\mathbb{P}(E)}(1))$. Then we have
$$\mathcal{L}=H^{\otimes\alpha}\otimes pr_2^*\mathcal{K},$$
for some line bundle $\mathcal{K}$ on $\mathbb{P}(E_\phi)$. Thus we
have
$$R^0pr_{2*}\mathcal{L}=Sym^\alpha(\pi^*E^*)
\otimes \mathcal{K}
=Sym^\alpha(\pi^*\mathcal{P}\oplus\pi^*\mathcal{P}^{-1}) \otimes
\mathcal{K},$$ where $\pi:\mathbb{P}(E_\phi)\rightarrow
T\times\widehat{T}$ is the structural map. Here $\alpha$ has to be
non negative, as $\mathcal{L}$ has a non zero section.

The non zero section of $\mathcal{L}$ defining $Z$ thus gives rise
to sections $\sigma_{\gamma\gamma'}$ of
$$\pi^*\mathcal{P}^{\otimes\gamma}\otimes \pi^*\mathcal{P}^{-\gamma'}\otimes\mathcal{K},$$
for $\gamma\geq0,\,\gamma'\geq0,\,\gamma+\gamma'=\alpha$.

Note that only one $\sigma_{\gamma\gamma'}$ can be non zero. Indeed,
by Lemma \ref{sub2}, the divisors of $\sigma_{\gamma\gamma'}$ have
to be combinations of $\Sigma_{1\phi}$ and $\Sigma_{2\phi}$ and the
two line bundles  $\mathcal{O}(\Sigma_{1\phi})$,
$\mathcal{O}(\Sigma_{2\phi})$ differ by a multiple of
$\pi^*\mathcal{P}_\phi$. Thus, if two sections
$\sigma_{\gamma\gamma'}$ were non zero, then we would  get a
proportionality relation between $\pi^*\mathcal{P}_\phi$ and
$\pi^*\mathcal{P}$ on $\mathbb{P}(E_\phi)$, which is not possible.

Thus there is only one non zero section $\sigma_{\gamma\gamma'}$.
There are now two possibilities: if the divisor $D_{\gamma\gamma'}$
of $\sigma_{\gamma\gamma'}$ is non-empty, then as $Z$ is irreducible
and contains $pr_2^{-1}D_{\gamma\gamma'}$, $Z$ must be a pull-back,
and
 Lemma
\ref{sub2} gives the result.

Next if the divisor $D_{\gamma\gamma'}$ of $\sigma_{\gamma\gamma'}$
is empty, one concludes that the line bundle $\mathcal{K}$ is a
pull-back :
$$\mathcal{K}=\pi^{*}\mathcal{K}'$$
for some line bundle $\mathcal{K}'$ on $T\times\widehat{T}$. But
then, $\mathcal{L}$ is also a pull-back:
$$\mathcal{L}=pr_1^{*}\mathcal{L}'$$
for some line bundle $\mathcal{L}'$ on $\mathbb{P}(E)$, and thus $Z$
is equal to $pr_1^{-1}(Z')$, for some $Z'\subset\mathbb{P}(E)$.
Lemma \ref{sub2} gives then the result.

The proof of b) is obtained by projecting codimension $2$ subsets of
$\mathbb{P}(E)\times_{T\times\widehat{T}}\mathbb{P}(E_\phi)$ to
$\mathbb{P}(E)$ and $\mathbb{P}(E_\phi)$.
\cqfd

The results above give us correspondingly the description of the codimension
$1$ and codimension $2$ analytic subsets of
$$Q:=\mathbb{P}(E)\times_{T\times\widehat{T}}\mathbb{P}(E_\phi)/
<(i,i_\phi),(\hat{i},\hat{i}_\phi)>$$
which dominate $K\times\widehat{K}$.

Namely they are the image in $Q$
 of the subvarieties described
above.

One interesting point is that the two hypersurfaces
$pr_1^{-1}\Sigma_1,\,pr_1^{-1}\Sigma_2$ descend to only one
irreducible hypersurface
\begin{eqnarray}
\label{Sigma}\Sigma\subset Q,
\end{eqnarray}
because the two factors in the splitting
$E=\mathcal{P}\oplus\mathcal{P}^{-1}$ are exchanged under $i$, so
that $ pr_1^{-1}\Sigma_1,\,pr_1^{-1}\Sigma_2$ are permuted by
$<(i,i_\phi),(\hat{i},\hat{i}_\phi)>$. For the same reason,
$pr_2^{-1}\Sigma_1^\phi,\,pr_2^{-1}\Sigma_2^\phi$ give rise to only
one hypersurface $\Sigma_\phi$.

Similarly the $4$ codimension $2$ subvarieties
$$pr_1^{-1}\Sigma_i\cap pr_2^{-1}\Sigma_j^\phi,\,i=1,\,2,\,j=1,\,2$$
descend to only one irreducible subvariety $W$ of $Q$, because they
are permuted by the group  $<(i,i_\phi),(\hat{i},\hat{i}_\phi)>$.

Thus $Q$, and hence $X$ contain only one irreducible codimension $2$
subvariety $W$ which dominates $K\times\widehat{K}$.

\vspace{0,5cm}

{\bf Proof of Proposition \ref{propess}.}
The proof is now immediate from the analysis above. Starting from $X$, the only
modifications which we can do, whose center dominates $K\times\widehat{K}$, is
to blow-up $W$, because in a quadric, there is no contractible curve.
In the blown-up variety, we have as divisors the exceptional divisors,
the proper transforms of the divisors $\Sigma,\,\Sigma_\phi$ and they are the only one.
Furthermore, the only codimension $2$ closed analytic subset dominating
$K\times\widehat{K}$ is the  union of two copies of $W$, indexed by the choice of
one of the divisors $\Sigma,\,\Sigma_\phi$, since
$W=\Sigma\cap\Sigma_\phi$. The same situation happens each time we blow-up
one copy of $W$ appearing in the previous step.

The key point is now the following:
If the map $\psi:X'\dasharrow X$ is not defined over the generic point of
$K\times \widehat{K}$, which we can see
as a birational map
between surface bundles over the generic point of $K\times \widehat{K}$, then
after a finite sequence of blow-ups of $X$ along codimension $2$
subsets dominating $K\times \widehat{K}$,
some divisor $D\subset \widetilde{X}$ in the blown-up variety must be generically
contractible over $K\times \widehat{K}$, that is
be made of a disjoint union of  rational curves of self-intersection $-1$ in the
generic surface $\widetilde{X}_t$, while this divisor $D$ projects to a divisor
in $X$. This follows from the factorization of birational map between surfaces
(see \cite{beauville}).

 But as this divisor dominates
 $K\times \widehat{K}$, it must be one of those described above,
that is  a proper transform of $\Sigma,\,\Sigma_\phi$. The
contradiction comes from the fact that
 after the blow-up of $W$, the
 proper transforms of $\Sigma$ and $\Sigma_\phi$ are families
 of rational curves of self-intersection $-2$, and this self-intersection can only decrease
 after further blow-ups. One the other hand, if we do not blow-up anything,
 these divisors are families of curves of self-intersection $0$, which do not contract.

\cqfd
\section{Proof of Theorem \ref{maincoh}\label{sec3}.}
In this section, we assume the hypotheses of Theorem \ref{maincoh}, namely,
$X'$ is bimeromorphically equivalent to $X$, and $Y$ is a compact K\"ahler
manifold such that there exists an isomorphism
$$\gamma:H^*(Y,\mathbb{Q})\cong H^*(X',\mathbb{Q})$$
of graded algebras. We want to prove that $Y$ cannot be projective.

Our argumentation will be based on the analysis of the algebra
$H^*(X',\mathbb{Q})$ made in the previous section, and on the following Lemma
\ref{delignelemma}
due to Deligne (see \cite{deligne}, \cite{voisin}, section 3) which was already heavily
used in the last section of \cite{voisin}.

Let $B^*$ be a finite dimensional graded algebra over $\mathbb{Q}$ and assume that
each $B^k$ carries a rational Hodge structure, compatible with the product, i.e.
the product map
$$B^k\otimes B^l\rightarrow B^{k+l}$$
is a morphism of Hodge structures.
Let $Z\subset B^k_\mathbb{C}$ be an algebraic subset defined by homogeneous equations
which can be formulated using only  the product structure on $B^*$. We have in mind, eg
$$Z=\{\alpha\in B^k_\mathbb{C},\,\alpha^2=0\}$$
or, which will be also used in the sequel, $Z'=Sing\,Z$, for $Z$ as above.
\begin{lemm}\label{delignelemma} (Deligne) For $Z$ as above, let
$Z_1\subset Z$ be an union of irreducible reduced components
of $Z$. Assume that the $\mathbb{C}$-vector space
$<Z_1>$ generated by $Z_1$ is defined over $\mathbb{Q}$, that
is $<Z_1>=Z_{1\mathbb{Q}}\otimes\mathbb{C}$, for some $\mathbb{Q}$-vector space
$Z_{1\mathbb{Q}}\subset B_\mathbb{Q}^k$. Then
$Z_{1\mathbb{Q}}$ is a rational sub-Hodge structure of $B_\mathbb{Q}^k$.
\end{lemm}

Our first step is the following (notations are as in the previous section):
\begin{prop}\label{subHodge1} Let $X'$, $Y$, $\gamma$ be as above.
Then $\gamma^{-1}(A^2_{1\mathbb{Q}})$ and
$\gamma^{-1}(A^2_{2\mathbb{Q}})$ are   rational sub-Hodge structures of
$H^2(Y,\mathbb{Q})$.
\end{prop}
{\bf Proof.} We give the proof for $\gamma^{-1}(A^2_{1\mathbb{Q}})$, the proof
for $\gamma^{-1}(A^2_{2\mathbb{Q}})$ is identical.

 We have only to explain how to recover
the space $A^2_{1\mathbb{C}}$ as generated by a certain algebraic subset
of $H^2(X',\mathbb{C})$ defined using only the algebra structure on $H^*(X',\mathbb{C})$,
since then, via $\gamma$,  we will then recover similarly
$\gamma^{-1}(A^2_{1\mathbb{C}})\subset H^2(Y,\mathbb{C})$ and then
by Deligne's Lemma \ref{delignelemma}, we will know that
$\gamma^{-1}(A^2_{1\mathbb{Q}})$ is a sub-Hodge structure of
$H^2(Y,\mathbb{Q})$.

We first use Proposition \ref{propcruciale}. It says that the irreducible
components of the algebraic subset
$$Z_1=\{\alpha_1+d,\,d\in D_\mathbb{C},\,\alpha_1\in A^2_{1\mathbb{C}},\,
\alpha_1^2=0,\,d^2=0,\,\alpha_1d=0\}$$
containing the
algebraic
subset
$$Z_{A^1}:=\{\alpha\in A^2_{1\mathbb{C}},\,
\alpha^2=0\},$$ are  irreducible components of
$$Z=\{\alpha\in H^2(X',\mathbb{C}),\,\alpha^2=0\}.$$
Next Lemma \ref{produitdiv} says us that if we denote by $D_1$  the
$\mathbb{Q}$-vector subspace
of $H^2(X',\mathbb{Q})$ defined as
$$D_1:=\{d\in D,\,d\alpha=0,\,\forall\alpha\in A^2_{1\mathbb{Q}}\},$$
the condition
$$\alpha_1d=0\,\,{\rm in}\,\, H^4(X',\mathbb{C}),$$
for some $$0\not=\alpha_1\in A^2_{1\mathbb{C}},\,d\in D_\mathbb{C},$$
 implies that
$d\in D_{1\mathbb{C}}:=D_{1}\otimes \mathbb{C}$.

Using this Lemma, we conclude that  the following algebraic subset of $H^2(X',\mathbb{C})$,
$$Z'_1=\{\alpha_1+d,\,d\in D_{1\mathbb{C}},\,\alpha_1\in A^2_{1\mathbb{C}},\,
\alpha_1^2=0,\,d^2=0\},$$
also satisfies the property that
its irreducible components containing $Z_{A_1}$ are irreducible components
of $Z$.
Note now that the vector space
$A^2_{1\mathbb{C}}$ is defined over $\mathbb{Q}$ and generated by its algebraic
subset
$Z_{A^1}$,
because $A^*_1$ is the exterior algebra $\bigwedge^{even}\Gamma^*_\mathbb{Q}$.

Thus, it remains only to show how to recover $Z_{A^1}$ from
$Z'_1$. This is done as follows. Let
$$D'_{1\mathbb{C}}\subset D_{1\mathbb{C}}$$
be the complex vector space generated by the  algebraic subset
$$Z_{D_1}:=\{d\in D_{1\mathbb{C}}, d^2=0 \}.$$
 $D'_{1\mathbb{C}}$ is defined over $\mathbb{Q}$,
 that is $$
 D'_{1\mathbb{C}}=D'_1\otimes\mathbb{C}$$
 for some rational subspace $D'_1\subset H^2(X',\mathbb{Q})$,
   because  $D_{1\mathbb{C}}$ is,
 and $Z_{D_1}$
is  defined over $\mathbb{Q}$.

If $D'_{1}=0$, there is nothing to say because then $Z'_1=Z_{A^1}$.
In general, the formula defining $Z'_1$ shows that it is the ``join'' of
$Z_{D_1}$ and $Z_{A_1}$ in $D'_1\oplus A^2_1$.

Assume first that $Z_{D_1}\not=D'_{1\mathbb{C}}$. In this case we
recover $Z_{A_1}$ as a component of the singular locus of $Z'_1$
because the join of two algebraic sets admits one of these algebraic
sets as an union of component of its singular locus unless the other
one is linear. So in this case, we recover $Z_{A_1}$ from the
algebra structure of $H^*(X',\mathbb{C})$ and this is finished.

It remains only to exclude the possibility that
\begin{eqnarray}\label{mauvais}D'_{1}\not=0,\,Z_{D_1}=D'_{1\mathbb{C}}.
\end{eqnarray}
This is done by the following argument :
assume (\ref{mauvais}) holds. As $D'_{1}$ is a $\mathbb{Q}$-vector space, there would be in
particular a non zero real element
$d\in D\subset H^{1,1}_\mathbb{R}(X')$ such that
$$d^2=0,\,d\alpha=0,\,\forall\alpha\in A^2_{1\mathbb{R}}.$$
But there exists also a non-zero $$\alpha\in A^{1,1}_{1\mathbb{R}}:
=H^{1,1}_\mathbb{R}(X')\cap A^2_{1\mathbb{R}}$$ such that
$\alpha^2=0$. It follows that the rank $2$ real vector space
$$B:=<d,\alpha>\subset H^{1,1}_\mathbb{R}(X')$$
satisfies the property:
 $$\forall \gamma\in B,\,\,\gamma^2=0.$$

But this contradicts the Hodge index theorem (cf \cite{voisinbook} I, 6.3.2)
because $X'$ is dominated by a K\"ahler compact manifold, and it follows
that for some element $c\in H^{4n}(X',\mathbb{R})$,
the intersection form
$$u\mapsto cu^2\in H^{4n+4}(X',\mathbb{R})=\mathbb{R}$$
has only one positive sign on $H^{1,1}_\mathbb{R}(X')$, hence cannot admit a rank
$2$ real isotropic subspace.
Thus (\ref{mauvais}) leads to a contradiction, and the proposition is proved.

\cqfd
\begin{coro}\label{coroD} With the same assumptions and notations, the subspace
$$\gamma^{-1}(D)\subset H^2(Y,\mathbb{Q})$$
is a rational sub-Hodge structure.
\end{coro}
{\bf Proof.}
 We use Lemma \ref{dercoh}, b), which says that
 $D\otimes\mathbb{C}$ is an irreducible component of the set
 $$Z'=\{d\in H^2(X',\mathbb{C}),\,d^3\beta=0,\,\forall \beta\in A^{4n-2}_\mathbb{C}\}.
 $$
 It follows that
 $\gamma^{-1}(D)\otimes\mathbb{C}$
 is an irreducible component of the set
 $$\gamma^{-1}(Z')=\{d\in H^2(Y,\mathbb{C}),\,d^3\beta=0,\,\forall \beta\in \gamma^{-1}(
 A^{4n-2}_\mathbb{C})\}.
 $$
 But we know as a consequence of Proposition
 \ref{subHodge1} that $\gamma^{-1}(
 A^{4n-2})$ is a rational sub-Hodge structure of $H^{4n-2}(Y,\mathbb{Q})$.
 Indeed, it is equal to the degree $4n-2$ piece of the
 subalgebra generated by
 $\gamma^{-1}(A^{2})$ and $\gamma^{-1}(A^{2})$
 is a rational sub-Hodge structure of $H^{2}(Y,\mathbb{Q})$.

  It follows that its annihilator
$$\gamma^{-1}(
 A^{4n-2})_0=\{\delta\in H^6(Y,\mathbb{Q}),\,\delta\beta=0,\,
 \forall \beta\in \gamma^{-1}(
 A^{4n-2}_\mathbb{C})\}$$
 is also a rational  sub-Hodge structure of $H^{4n-2}(Y,\mathbb{Q})$.

 Hence there is an induced rational Hodge structure on the quotient
 $$H^6(Y,\mathbb{Q})/\gamma^{-1}(A^{4n-2})_0$$
 and we can apply Deligne's Lemma \ref{delignelemma} to
 the product
 $$H^2(Y,\mathbb{Q})^{\otimes3}\rightarrow H^6(Y,\mathbb{Q})/\gamma^{-1}(A^{4n-2})_0,$$
 which is compatible with the induced Hodge structure: Indeed,
 for this product, we have that
$\gamma^{-1}(D)\otimes\mathbb{C}$
 is an irreducible component of the set
 $$Z''=\{\delta\in H^2(Y,\mathbb{C}),\,\delta^3=0\}.$$
 As $\gamma^{-1}(D)$ is a rational subspace of $H^2(Y,\mathbb{Q})$, Lemma \ref{delignelemma}
 says that
 it is a rational sub-Hodge structure of $H^2(Y,\mathbb{Q})$.
\cqfd
{\bf Proof of Theorem \ref{maincoh}.}
The isomorphism of graded algebras
$$\gamma: H^*(Y,\mathbb{Q})\cong H^*(X',\mathbb{Q})$$
must be compatible up to a coefficient with Poincar\'e duality, which is given by the cup-product
and  isomorphisms
$$ H^{4n+4}(X',\mathbb{Q})=\mathbb{Q},\,\,H^{4n+4}(Y,\mathbb{Q})=\mathbb{Q}.$$
As $\gamma^{-1}(A^2)$ is a rational sub-Hodge structure of
$H^2(Y,\mathbb{Q})$, so is
\begin{eqnarray}\label{sansnom}
\gamma^{-1}(A^{4n-4})\subset
H^{4n-4}(Y,\mathbb{Q}),
\end{eqnarray}
because it is equal to the degree $4n-4$ piece of the subalgebra of
$H^{*}(Y,\mathbb{Q})$ generated by $\gamma^{-1}(A^2)$.

Now, the map which is Poincar\'e dual to the inclusion
$$A^{4n-4}\subset
H^{4n-4}(X',\mathbb{Q})$$
is the map
$$(q\circ\psi)_*:H^8(X',\mathbb{Q})\rightarrow H^4((T/<\pm Id>)\times
(\widehat{T}/<\pm Id>),\mathbb{Q})$$
$${\cong} A_{1\mathbb{Q}}^4\oplus A_{1\mathbb{Q}}^2
\otimes A_{2\mathbb{Q}}^2\oplus A_{2\mathbb{Q}}^4,$$
where the last isomorphism is given by the K\"unneth decomposition.
We shall denote by
$$\kappa:H^4((T/<\pm Id>\times(\widehat{T}/<\pm Id>),\mathbb{Q})\rightarrow
 A_{1\mathbb{Q}}^2
\otimes A_{2\mathbb{Q}}^2$$
the K\"unneth projector given by the decomposition above.

Applying $\gamma^{-1}$, we thus get a projection
$$ H^8(Y,\mathbb{Q})\rightarrow
\gamma^{-1}(A_{1\mathbb{Q}}^4)\oplus \gamma^{-1}(A_{1\mathbb{Q}}^2)
\otimes \gamma^{-1}(A_{2\mathbb{Q}}^2)\oplus \gamma^{-1}(A_{2\mathbb{Q}}^4)$$
which must be a morphism of Hodge structures as its transpose (\ref{sansnom})
is.
Composing further with the projection (conjugate via $\gamma$ to $\kappa$)
$$\gamma^{-1}(A_{1\mathbb{Q}}^4)\oplus \gamma^{-1}(A_{1\mathbb{Q}}^2)
\otimes \gamma^{-1}(A_{2\mathbb{Q}}^2)\oplus \gamma^{-1}(A_{2\mathbb{Q}}^4)
\rightarrow \gamma^{-1}(A_{1\mathbb{Q}}^2)
\otimes \gamma^{-1}(A_{2\mathbb{Q}}^2),$$
which is also a morphism of Hodge structures because
$\gamma^{-1}(A_{1\mathbb{Q}}^2)$ and
$ \gamma^{-1}(A_{2\mathbb{Q}}^2)$ are sub-Hodge structures of $H^2(Y,\mathbb{Q})$,
we get finally
a morphism of Hodge structures
$$H^8(Y,\mathbb{Q})\rightarrow \gamma^{-1}(A_{1\mathbb{Q}}^2)
\otimes \gamma^{-1}(A_{2\mathbb{Q}}^2).$$
Restricting it to the sub-Hodge structure $\gamma^{-1}(D)^4=\gamma^{-1}(D^4)\subset H^8(Y,\mathbb{Q})$
generated by $\gamma^{-1}(D)$, we finally get a morphism of rational Hodge structures
$$\pi_\gamma:\gamma^{-1}(D^4)\rightarrow
\gamma^{-1}(A_{1\mathbb{Q}}^2)
\otimes \gamma^{-1}(A_{2\mathbb{Q}}^2),$$
which is conjugate  via $\gamma$ to
the restriction of $\kappa\circ(q\circ\psi)_*$ to $D^4$.

We have now the following two Lemmas :
\begin{lemm}\label{avantder} The image of
$$\kappa\circ(q\circ\psi)_*: D^4\rightarrow
A_{1\mathbb{Q}}^2
\otimes A_{2\mathbb{Q}}^2$$
contains
$$Id\in Hom\,(A_{1\mathbb{Q}}^2,A_{1\mathbb{Q}}^2)\cong A_{1\mathbb{Q}}^2
\otimes A_{2\mathbb{Q}}^2$$
and
$$\phi^*=\wedge^2{^t\phi}\in Hom\,(A_{1\mathbb{Q}}^2,A_{1\mathbb{Q}}^2)\cong A_{1\mathbb{Q}}^2
\otimes A_{2\mathbb{Q}}^2.$$
\end{lemm}
Let now
$$\Pi_\gamma=\gamma^{-1}\otimes\gamma^{-1}(Im\,\kappa\circ(q\circ\psi)_*)\subset
\gamma^{-1}(A_{1\mathbb{Q}}^2)\otimes\gamma^{-1}(A_{2\mathbb{Q}}^2)$$
be
the image of $\pi_\gamma$.
\begin{lemm}\label{derlemma}   a) The generic element of $\Pi_\gamma$ is non-degenerate.

(Here we see $u\in \gamma^{-1}(A_{1\mathbb{Q}}^2)\otimes\gamma^{-1}(A_{2\mathbb{Q}}^2)$
as an element of
$$Hom\,(\gamma^{-1}(A_{2\mathbb{Q}}^{2})^*,\gamma^{-1}(A_{1\mathbb{Q}}^2))$$
and non-degenerate means invertible.)

b)
The $\mathbb{Q}$ vector subspace $\Pi_\gamma'$ of
$End\,(\gamma^{-1}(A_{2\mathbb{Q}}^{2*}) )$ generated by
the
$u^{-1}\otimes v,\,u$ non-degenerate in $\Pi_\gamma$, consists of Hodge classes
in $End\,(\gamma^{-1}(A_{1\mathbb{Q}}^{2*})
)$, (relative to the Hodge structures
on $\gamma^{-1}(A_{2\mathbb{Q}}^{2})$
induced  by the Hodge structure on $H^2(Y,\mathbb{Q})$).
\end{lemm}

Assuming these Lemmas, the proof is now concluded as follows.

The two Lemmas together imply that
 the Hodge structure on $\gamma^{-1}(A_{2\mathbb{Q}}^{2*})$
admits an endomorphism conjugate to $ \phi_T^*=\wedge^2{^t\phi}$.
Hence dually the Hodge structure on
$\gamma^{-1}(A_{2\mathbb{Q}}^{2})$ admits a morphism conjugate to
$\wedge^2{\phi}$.

The proof concludes then exactly as in \cite{voisin}, 3.2: The above implies that either
the Hodge structure on
$\gamma^{-1}(A_{2\mathbb{Q}}^2)$ is trivial or it does not contain any Hodge class.
The first case is excluded by a Hodge index argument.

Next, working symmetrically with $A_{1\mathbb{Q}}^{2}$, we conclude similarly that
the Hodge structure on
$\gamma^{-1}(A_{2\mathbb{Q}}^2)$ does not contain any Hodge class.

Thus it follows from
Corollary \ref{coroD} that the only degree $2$ Hodge classes on $Y$
are contained in $\gamma^{-1}(D)$.

But  we look now at the intersection
form
$$q_d=\int_Y d^{4n}\alpha\beta$$
for $d\in \gamma^{-1}(D)$, and  we conclude that it is zero on
$\gamma^{-1}(A^2_{1\mathbb{Q}})$, because the same is true for $D$
and $A^2_{1\mathbb{Q}}$ on $X'$. Thus for no degree $2$ Hodge class
$d$ on $Y$, the sub-Hodge structure
$\gamma^{-1}(A_{1\mathbb{Q}}^2)\subset H^2(Y,\mathbb{Q})$ can be
polarized by $q_d$. Thus by \cite{voisinbook}, I, 6.3.2, $Y$ cannot
be projective.

\cqfd
{\bf Proof of Lemma \ref{avantder}.}
We first reduce to the case where $X'=X$:First of all, using
Lemma \ref{sub1}, we conclude that for any non-empty
Zariski open set $U$ of $K_0\times\widehat{K}_0$,
the restriction map
$$ H^4(K_0\times\widehat{K}_0,\mathbb{Q})=
H^4((T/\pm Id)\times(\widehat{T}/\pm Id),\mathbb{Q}) \rightarrow
H^4(U,\mathbb{Q})$$ is an isomorphism. Now we have the commutative
diagram:
\begin{eqnarray}\label{restUmatrix}
\begin{matrix}&D^4\subset H^8(X',\mathbb{Q})&\stackrel{rest_U}{\rightarrow}&
H^8(X'_U,\mathbb{Q}) \\
&(q\circ\psi)_*\downarrow&&\downarrow(q\circ\psi)_*^U\\
&H^4((T/\pm Id)\times(\widehat{T}/\pm Id),\mathbb{Q})&
\stackrel{rest_U}{\cong}& H^4(U,\mathbb{Q}).
\end{matrix}
\end{eqnarray}

 We use now Proposition \ref{propess} which says that the meromorphic map
 $\psi$ is well defined on a Zariski open set $X'_U$ as above.
 We thus have a commutative diagram:
$$\begin{matrix}&D^4_{X\mid X_U}&\stackrel{\psi_U^*}{\rightarrow}&
D^4_{\mid X'_U} \\
&q_{U*}\downarrow\,\,\,\,\,\,&&\,\,\,\,\,\,\,\,\,\,\,\,\,\,\,\,\,\,\,\,\downarrow(q\circ\psi)_{U*}\\
&H^4(U,\mathbb{Q})&
\cong& H^4(U,\mathbb{Q}),
\end{matrix}
$$
where $q_U,\,(q\circ\psi)_U$ denote the restrictions of
$q,\,q\circ\psi$ to $X_U,\,X'_U$ respectively. We used here the fact
that degree $2$ Hodge classes on $X$, restricted to $X_U$, pull-back
via $\psi_U$ to degree $2$ Hodge classes on $X'$, restricted to
$X'_U$, which follows from the fact that $\psi$ is meromorphic.

Writing for $X$ the same diagram as (\ref{restUmatrix}),
we conclude that it suffices to prove the result for
$X$.

Next, we look at the following Cartesian diagram:
$$
\begin{matrix}&\tilde{q}:&\mathbb{P}(E)_0\times_{T_0\times\widehat{T_0}}\mathbb{P}(E_\phi)_0&
\rightarrow& T_0\times\widehat{T_0}\\
&&e\downarrow&&e\downarrow\\
&q:&X_0&\rightarrow& K_0\times\widehat{K_0}
\end{matrix},
$$
where the lower indices $0$ denote the restrictions of the
projective bundles to $T_0\times\widehat{T_0}$, the vertical maps
denoted by $e$ are the quotient maps, and the induced map
$$H^4(K_0\times\widehat{K_0},\mathbb{Q})\rightarrow H^4(T_0\times\widehat{T_0},\mathbb{Q})$$
are injective. Here $X_0$ is the Zariski open set of $X$ which is
the smooth part of the quotient $Q$. Arguing as before, we see that
we can replace $X$ by $X_0$, and then $X_0$ by its \'etale cover
$\mathbb{P}(E)_0\times_{T_0\times\widehat{T_0}}\mathbb{P}(E_\phi)_0$.
 Thus the
result for $X$ follows from the following formulas (\ref{formule}):

Let $\Sigma,\,\Sigma_\phi$ be the two divisors of
(\ref{Sigma}), and let $s,\,s_\phi\in Hdg^2(X,\mathbb{Q})$ be their cohomology classes.
Then we have
\begin{eqnarray}
\label{formule}
\tilde{q}_*(e^*(s^3s_\phi))= 16 Id\in Hom\,(H^2(T_0,\mathbb{Q}),H^2(T_0,\mathbb{Q}))=
H^2(T_0,\mathbb{Q})\otimes H^2(\widehat{T}_0,\mathbb{Q}),\\ \nonumber
\tilde{q}_*(e^*(s s_\phi^3))= 16 \phi^*\in Hom\,(H^2(T_0,\mathbb{Q}),H^2(T_0,\mathbb{Q}))=
H^2(T_0,\mathbb{Q})\otimes H^2(\widehat{T}_0,\mathbb{Q}).
\end{eqnarray}
This is computed as follows: let $s_1,\,s_2$ be the classes of
the divisors $\Sigma_1,\,\Sigma_2$ of $\mathbb{P}(E)\times_{T_0\times\widehat{T_0}}\mathbb{P}(E_\phi)$
given by the
decomposition $E=\mathcal{P}\oplus\mathcal{P}^{-1}$ and similarly
let
$s_1^\phi,\,s_2^\phi$ be the classes of
the divisors $\Sigma_1^\phi,\,\Sigma_2^\phi$
 of $\mathbb{P}(E)\times_{T_0\times\widehat{T_0}}\mathbb{P}(E_\phi)$
given by the
decomposition $E_\phi=\mathcal{P}_\phi\oplus\mathcal{P}^{-1}_\phi$.
Then we have
$$e^*(s)=s_1+s_2,\,e^*(s_\phi)=s_1^\phi+s_2^\phi.$$
Let $h$, $h_\phi$ be respectively
$c_1(\mathcal{O}_{\mathbb{P}(E)}(1))$,
$c_1(\mathcal{O}_{\mathbb{P}(E_\phi)}(1))$, or rather their
pull-backs to the fibered product
$\mathbb{P}(E)\times_{T_0\times\widehat{T_0}}\mathbb{P}(E_\phi)$.
Let $p,\,p_\phi$ be the  classes
$c_1(\mathcal{P}),\,c_1(\mathcal{P}_\phi)$. Then we have
$$s_1=\tilde{q}^*p-h,\,s_2=-\tilde{q}^*p-h,$$
$$s_1^\phi=\tilde{q}^*p_\phi-h_\phi,\,s_2^\phi=-\tilde{q}^*p_\phi-h_\phi.$$
Thus
$$e^*(s)=-2h,\,e^*(s_\phi)=-2h_\phi,$$
and
$$e^*(s^3s_\phi)=16h^3h_\phi,\,e^*(s^3_\phi s)=16h^3_\phi h.$$
Applying $\tilde{q}_*$ we conclude that
$$\tilde{q}_*(e^*(s^3s_\phi))=-16c_2(E),\,
\tilde{q}_*(e^*(s^3_\phi s))=-16c_2(E_\phi).$$
As $E=\mathcal{P}\oplus\mathcal{P}^{-1}$, and
$E_\phi=\mathcal{P}_\phi\oplus\mathcal{P}_\phi^{-1}$, it follows that
$$c_2(E)=-p^2,\,c_2(E_\phi)=-p_\phi^2.$$
Finally, since $p$ identifies to
$$Id\in Hom\,(H^1(T,\mathbb{Q}),H^1(T,\mathbb{Q}))=H^1(T,\mathbb{Q})
\otimes H^1(\widehat{T},\mathbb{Q})\subset H^2(T\times\widehat{T},\mathbb{Q}),$$
we get that $p^2$
identifies to
$$Id\in Hom\,(H^2(T,\mathbb{Q}),H^2(T,\mathbb{Q}))=H^2(T,\mathbb{Q})
\otimes H^2(\widehat{T},\mathbb{Q})\subset H^4(T\times\widehat{T},\mathbb{Q}),$$
and similarly
$p^2_\phi$
identifies to
$$Id\in Hom\,(H^2(T,\mathbb{Q}),H^2(T,\mathbb{Q}))=H^2(T,\mathbb{Q})
\otimes H^2(\widehat{T},\mathbb{Q})\subset H^4(T\times\widehat{T},\mathbb{Q}).$$
Thus (\ref{formule}) is proved, which concludes the proof of the Lemma.
\cqfd
{\bf Proof of Lemma \ref{derlemma}.} The first statement is obvious
by Lemma \ref{avantder}.

Next, because we proved that $\Pi_\gamma$ is a sub-Hodge structure of
$$\gamma^{-1}(A_{1\mathbb{Q}}^2)\otimes\gamma^{-1}(A_{2\mathbb{Q}}^2),$$
 it follows that
the space $\Pi'_\gamma$ is a sub-Hodge structure of
$End\,(\gamma^{-1}(A_{2\mathbb{Q}}^{2*}) )$, and thus, so is the sub-algebra of
$End\,(\gamma^{-1}(A_{2\mathbb{Q}}^{2*}) )$ generated by $\Pi'_\gamma$.
On the other hand, $\Pi'_\gamma$ is conjugate
via $^t\gamma$ to the corresponding subspace of
$End\,(A_{2\mathbb{Q}}^{2*})$, defined similarly starting from
$Im\,\kappa\circ(q\circ\psi)_{*\mid D^4}$. This last subspace is contained in the
space of endomorphisms of Hodge structures of $A_{2\mathbb{Q}}^{2*}$,
which has been computed to be equal to the algebra generated
by $\phi_{\widehat{T}*}=\wedge^2\phi$ (see proof of Lemma \ref{plustard2}).

The key point is that because $\wedge^2\phi$ is diagonalizable, this
algebra tensored with $\mathbb{C}$ has no nilpotent element. It
follows that $\Pi'_\gamma\otimes \mathbb{C}$ has no nilpotent
element. But as $\Pi'_\gamma$ is a sub-Hodge structure of
$End\,(\gamma^{-1}(A_{2\mathbb{Q}}^{2*}))$, it follows that it is
pure of type $(0,0)$, that is made of Hodge classes, because
elements of type $(-k,k),\,k>0$ are nilpotent.

\cqfd

\end{document}